# Box-counting dimension of a kind of fractal interpolation surface on rectangular grids

CholHui Yun and MunChol Kim
Faculty of Mathematics, Kim Il Sung University, D. P. R. Korea
Email: yuncholhui@yahoo.com

**Abstract**

We estimate a Box-counting dimension of fractal surfaces which are generated by iterated function systems with a vertical contraction factor function on an arbitrary data set over rectangular grids and can express well a lot of natural surfaces with very complicated structures.

**Keywords:** Iterated function system (IFS); Fractal interpolation function(FIF); Fractal interpolation surface(FIS); Box-counting dimension

## 1. Introduction

A fractal interpolation surface (FIS) is a fractal set which is a graph of an interpolation function. Therefore, constructions of FISs are closely associated with one of fractal interpolation function, i.e. interpolation functions whose graph is fractal sets.

By Barnsley[3], FIFs were introduced in 1986 and after that have widely been studied in approximation theory, image compression, computer graphics and so on..

In many papers constructions of fractal interpolation surfaces on the basis of IFSs were studied. (See [5, 6, 7, 11, 13, 14, 17]) Massopust presented the construction of FISs on rectangular data sets, at which the interpolation points on the boundary are coplanar. Geronimo and Hardin generalized this construction to allow more general boundary data.

Constructions of fractal interpolation surfaces, which interpolate a data set over rectangular grid, were studied in [5, 6, 11, 14]. A lack of constructions is to use IFSs constructed with a restricted data set, whose data points on the boundary are collinear, constant contraction factor, quadratic polynomials. In [14], these constructions were generalized with an arbitrary data set, Lipschitz function, contraction factor function, lower and upper bounds for the Box-counting of the constructed surface.

In this paper, we consider the problem on improving estimation of Box counting dimension of fractal surfaces with vertical contraction factor function as fractal surfaces constructed in [14].

## 2. FISs over the rectangular grids

We consider the construction of fractal surfaces presented in [14]. Let the data set over the rectangular grid be

$$P = \{(x_i, y_j, z_{ij}) \in \mathrm{P}^3; i = 0, 1, \cdots, n,\ j = 0, 1, \cdots, m\}$$
$$(x_0 < x_1 < \cdots < x_n,\ y_0 < y_1 < \cdots < y_m)$$

and denote

$$N_{nm} = \{1, \cdots, n\} \times \{1, \cdots, m\},\ I_x = [x_0, x_n],\ I_y = [y_0, y_m], E = I_x \times I_y,$$
$$I_{x_i} = [x_{i-1}, x_i],\ I_{y_j} = [y_{j-1}, y_j],\ E_{ij} = I_{x_i} \times I_{y_j}, (i, j) \in N_{nm},$$
$$P_{x_\alpha} = \{(x_\alpha, y_l, z_{\alpha l}) \in P; l = 0, 1, \cdots, m\},\quad \alpha = 0, 1, \cdots, n,$$
$$P_{y_\beta} = \{(x_k, y_\beta, z_{k\beta}) \in P; k = 0, 1, \cdots, n\},\ \beta = 0, 1, \cdots, m.$$

We define the domain contraction transformations $L_{ij} : E \to E_{ij}$ for $(i, j) \in N_{nm}$ by

$$L_{ij}(x, y) = (L_{x_i}(x), L_{y_j}(y))$$

where $L_{x_i} : I_x \to I_{x_i}$ and $L_{y_j} : I_y \to I_{y_j}$ are contractive homeomorphisms with contraction factors $c_{x_i}, c_{y_j}$ satisfying the following conditions.

(1) $L_{x_i} : \{x_0, x_n\} \to \{x_{i-1}, x_i\},\ L_{y_j} : \{y_0, y_m\} \to \{y_{j-1}, y_j\}$

(2) For any $i (\in \{1, \cdots, n-1\}), j (\in \{1, \cdots, m-1\})$, there are $x_k (\in \{x_0, x_n\})$, $y_l (\in \{y_0, y_m\})$ such that

$$L_{x_{i+1}}(x_k) = L_{x_i}(x_k) = x_i,\quad L_{y_{j+1}}(y_l) = L_{y_j}(y_l) = y_j. \tag{1}$$

Then $c_{ij} = \mathrm{Max}\{c_{x_i}, c_{y_j}\}, (i, j) \in N_{nm}$ are contractivity factors of the transformations $L_{ij}$.

Functions $F_{ij}: E\times \mathbb{R} \to \mathbb{R}, (i, j)\in N_{nm}$ are defined by

$$F_{ij}(x, y, z) = s(L_{ij}(x, y))(z - g(x, y)) + h(L_{ij}(x, y)),$$

where $s(x, y)$ is a vertical continuous contraction function such that $0<|s(x, y)|<1$ on $E$, $h(x, y)$, $g(x, y)$ are continuous Lipschitz functions on $E$ with the Lipschitz constants $L_h, L_g$ satisfying

$$g(x_\alpha, y_\beta) = z_{\alpha\beta},\ (\alpha, \beta)\in\{0, n\}\times\{0, m\},$$
$$h(x_i, y_j) = z_{ij},\ (i, j)\in\{0,1,\cdots,n\}\times\{0,1,\cdots,m\}.$$

Then transformations $W_{ij} = (L_{ij}, F_{ij})^T, (i, j)\in N_{nm}$ coincide on common borders and are contractions with respected to some metric that is equivalent to the Euclidean metric on $\mathbb{R}^3$. And the attractor A of the IFS $\{\mathbb{R}^3; W_{ij}, i=1,\cdots,n,\ j=1,\cdots,m\}$ is a graph of a continuous function $f: E\to\mathbb{R}$, i.e. a surface in $\mathbb{R}^3$.

The type of $f$ is as follows;

$$\begin{aligned} f(x, y) &= s(x, y) f(L_{ij}^{-1}(x, y)) + Q(x, y) \\ Q(x, y) &= -s(x, y) g(L_{ij}^{-1}(x, y)) + h(x, y) \end{aligned}, \tag{2}$$

## 2. Box-counting dimension of interpolation surface

In this section, we calculate lower and upper bounds for Box-counting dimension of the surface constructed as above with the data set

$$P = \left\{\left(x_0 + \frac{x_n - x_0}{n} i, y_0 + \frac{y_n - y_0}{n} j, z_{ij}\right)\in\mathbb{R}^3; i, j = 0,1,\cdots,n\right\}$$

Since there is a bi-Lipschitz mapping between some rectangular in $\mathbb{R}^2$ and $[0,1]\times[0,1]$, and Box-counting dimension is invariant under bi-Lipschitz mapping, we can assume that $E=[0,1]\times[0,1]$. Then $P = \left\{\left(\frac{i}{n}, \frac{j}{n}, z_{ij}\right)\in\mathbb{R}^3; i, j = 0,1,\cdots,n\right\}$.

For $D(\subset\mathbb{R}^2)$, let the maximum range of $f$ on $D$ be denoted by

$R_f[D] := \sup\{|f(x_2, y_2) - f(x_1, y_1)|; (x_1, y_1), (x_2, y_2)\in D\}$

**[Lemma]** Let $W: D\times\mathbb{R}\to D\times\mathbb{R}$ be the form

$$W\begin{pmatrix} x \\ y \\ z \end{pmatrix} = \begin{pmatrix} L(x, y) \\ F(x, y, z) \end{pmatrix} = \begin{pmatrix} L_x(x) \\ L_y(y) \\ s(L(x, y))z + Q(x, y) \end{pmatrix}$$

where $Q$ is Lipschitz function with the Lipschitz constants $L_Q$ , $L$ is the domain contraction transformation with contraction factor $c_L$ defined as $L_{ij}$ and $s(x, y)$ is also contraction function with $0 < |s(x, y)| < 1$ .Then, for any continuous function $f : D \to \mathrm{P}$ ,

$$R_{F(L^{-1},\, f \circ L^{-1})}[L(D)] \le \bar{s} R_f[D] + \mathrm{diam}(D)(c_s \overline{f} + L_Q) .$$

Where $\mathrm{diam}(D)$ is a diameter of the set $D$, $\bar{s} = \underset{D}{\mathrm{Max}} |s(x, y)|$, $c_s$ is a contraction factor of $s(x, y)$, $\overline{f} = \underset{D}{\mathrm{Max}} |f(x, y)|$.

(Proof) For $(x, y), (x', y') (\in L(D))$, let denote $L^{-1}(x, y) =: (\tilde{x}, \tilde{y})$, $L^{-1}(x', y') =: (\tilde{x}', \tilde{y}') (\in D)$. Then

$$\begin{aligned}
&| F(L^{-1}, f \circ L^{-1})(x, y) - F(L^{-1}, f \circ L^{-1})(x', y') | = \\
&= | F(L^{-1}(x, y), f \circ L^{-1}(x, y)) - F(L^{-1}(x', y'), f \circ L^{-1}(x', y')) | \\
&= | s(x, y) f(\tilde{x}, \tilde{y}) + Q(\tilde{x}, \tilde{y}) - s(x', y') f(\tilde{x}', \tilde{y}') - Q(\tilde{x}', \tilde{y}') | \\
&= | s(x, y) f(\tilde{x}, \tilde{y}) - s(x, y) f(\tilde{x}', \tilde{y}') + s(x, y) f(\tilde{x}', \tilde{y}') - s(x', y') f(\tilde{x}', \tilde{y}') \\
&+ Q(\tilde{x}, \tilde{y}) - Q(\tilde{x}', \tilde{y}') | \\
&\le \bar{s} R_f[D] + c_s d((x, y), (x', y')) \overline{f} + L_Q d((\tilde{x}, \tilde{y}), (\tilde{x}', \tilde{y}')) \\
&\le \bar{s} R_f[D] + \mathrm{diam}(D)(c_s \overline{f} + L_Q).
\end{aligned}$$

(end of proof)

For $N \times N$ matrix $U = (u_{ij})_{N \times N}$, $V = (v_{ij})_{N \times N}$, a relation "<" is defined by

$$U < V \overset{d}{\Leftrightarrow} u_{ij} < v_{ij},\ i, j = 1, 2, \cdots, N .$$

And let denote $\bar{s}_{ij} = \underset{E_{ij}}{\mathrm{Max}} |s(x, y)|$, $\tilde{s}_{ij} = \underset{E_{ij}}{\mathrm{Min}} |s(x, y)|$.

**[Theorem]** If there is $\alpha_0 (\in \{0, 1, \cdots, n\})$ (or $\beta_0 (\in \{0, 1, \cdots, n\})$) such that the points of $P_{x_{\alpha_0}}$ (or $P_{y_{\beta_0}}$) are not collinear, then the Box-dimension $\dim_B A$ of the

attractor $A$ is as follows;

(1) If $\sum_{i,j=1}^{n} \tilde{s}_{ij} > n$, then

$$1+\log_n^{\sum_{i,j=1}^{n} \tilde{s}_{ij}} \le \dim_B A \le 1+\log_n^{\sum_{i,j=1}^{n} \bar{s}_{ij}} .$$

(2) If $\sum_{i,j=1}^{n} \bar{s}_{ij} \le n$,

$$\dim_B A = 2 .$$

(Proof) Proof of (1). By the hypothesis of the theorem, maximum vertical distance calculated only with respected to Z axis from the points of $P_{x_{\alpha_0}}$ (or $P_{y_{\beta_0}}$ ) to the line through two end points of $P_{x_{\alpha_0}}$ (or $P_{y_{\beta_0}}$ ) is positive. This is called a height and denoted by $H$ .

After each $W_{ij}$ is applied to the interpolation points in $E$, we obtain $(n+1)^2$ new points in every $E_{ij}$ and the vertical lines parallel to $z$ axis are mapped to the vertical lines parallel to $z$ axis by $W_{ij}$. Hence all of vertical lines with length $H$ are mapped to vertical lines in $E_{ij}$ whose length is more than $\tilde{s}_{ij}H$ .

And using lemma for $R_f[E_{ij}]$, we have

$$R_f[E_{ij}] \le \bar{s}_{ij} R_f[E] + b,$$

where $b=\sqrt{2}(c_s \overline{f} + L_Q)$. Let denote $R_f[E_{ij}]$ by $R_{ij}$.

Let a injection $\tau : \{1,\cdots,n\}\times\{1,\cdots,n\} \to \{1,\cdots,n^2\}$ be define by $\tau(i,\ j)=(i-1)n+j$. Let $n^2 \times n^2$ diagonal matrix $\bar{S},\ \tilde{S}$ and vectors $\mathbf{h}_1, \mathbf{r}, \mathbf{u}_1, \mathbf{i}$ be as follows;

$$\bar{S} = \mathrm{diag}\left(\bar{s}_{\tau^{-1}(1)},\cdots,\bar{s}_{\tau^{-1}(n^2)}\right),\quad \tilde{S} = \mathrm{diag}\left(\tilde{s}_{\tau^{-1}(1)},\cdots,\tilde{s}_{\tau^{-1}(n^2)}\right),$$

$$\mathbf{h}_1 = \begin{pmatrix} \tilde{s}_{\tau^{-1}(1)}H \\ \vdots \\ \tilde{s}_{\tau^{-1}(n^2)}H \end{pmatrix},\quad \mathbf{r} = \begin{pmatrix} \bar{s}_{\tau^{-1}(1)}R_{\tau^{-1}(1)} \\ \vdots \\ \bar{s}_{\tau^{-1}(n^2)}R_{\tau^{-1}(n^2)} \end{pmatrix},\quad \mathbf{i} = \begin{pmatrix} 1 \\ \vdots \\ 1 \end{pmatrix},\quad \mathbf{u}_1 = \mathbf{r} + b\mathbf{i} .$$

For $r(>0)$ let denote $\varepsilon_r := \dfrac{1}{n^r}$. Then $\varepsilon_r \to 0 \Leftrightarrow r \to \infty$. And let $N(\varepsilon_r)$ be defined by the smallest number of $\varepsilon_r$-mesh cubs that cover $A$. Then, since $A$ is the graph of continuous function on $E$, the smallest number of $\varepsilon_r$-mesh cubs that cover $(E_{ij} \times \mathrm{P}) \bigcap A$ is greater than one of $\varepsilon_r$-mesh cubs that cover vertical lines with the length $\tilde{s}_{ij} H$ ,and is less than one of $\varepsilon_r$-mesh cubs that cover parallelepiped $E_{ij} \times [\tilde{f}_{ij}, \overline{f}_{ij}]$. Where $\tilde{f}_{ij} = \operatorname*{Min}_{E_{ij}} | f(x, y) |$, $\overline{f}_{ij} = \operatorname*{Max}_{E_{ij}} | f(x, y) |$. Hence

$$\sum_{i,j=1}^{n} [\tilde{s}_{ij} H \varepsilon_r^{-1}] \le N(\varepsilon_r) \le \sum_{i,j=1}^{n} ([(\overline{s}_{ij} R_{ij} + b)\varepsilon_r^{-1}] + 1) \left( \left[ \frac{\varepsilon_r^{-1}}{n} \right] + 1 \right)^2 ,$$

$$\sum_{k=1}^{n^2} (\tilde{s}_{\tau^{-1}(k)} H \varepsilon_r^{-1}) - n^2 \le N(\varepsilon_r) \le \sum_{k=1}^{n^2} ((\overline{s}_{\tau^{-1}(k)} R_{\tau^{-1}(k)} + b)\varepsilon_r^{-1} + 1) \left( \left[ \frac{\varepsilon_r^{-1}}{n} \right] + 1 \right)^2$$

that is,

$$\Phi(\mathbf{h}_1 \varepsilon_r^{-1}) - n^2 \le N(\varepsilon_r) \le \Phi(\mathbf{u}_1 \varepsilon_r^{-1} + \mathbf{i}) \left( \left[ \frac{\varepsilon_r^{-1}}{n} \right] + 1 \right)^2 ,$$

where for vector $\mathbf{a} = (a_1, \cdots, a_m)$, $\Phi(\mathbf{a}) = a_1 + \cdots + a_m$ and $\varepsilon_r^{-1} \ge n$.

After applying $W_{ij}$ to $E$ two times, we get $n^2$ squares of side $1/n^2$ in $E_{ij}$. Since each square is obtained from each $E_{kl}$ lying inside $E$ by transformations $W_{kl}$, the sum of maximum ranges of $f$ on $n^2$ squares of side $1/n^2$ contained in $E_{ij}$ is less than or equal to coordinates of vector

$$\mathbf{u_2} = \overline{S} C \mathbf{u_1} + nb\mathbf{i}$$

and greater than or equal to coordinates of vector

$$\mathbf{h_2} = \tilde{S} C \mathbf{h_1}$$

And we have

$$\Phi(\mathbf{h_2}\varepsilon_r^{-1}) - n^4 \le N(\varepsilon_r) \le \Phi(\mathbf{u_2}\varepsilon_r^{-1} + n^2\mathbf{i})\left(\left[\frac{\varepsilon_r^{-1}}{n^2}\right]+1\right)^2$$

,where $C$ is $n^2 \times n^2$ matrix whose entries are all 1 and $\varepsilon_r^{-1} \ge n^2$.

After applying $W_{ij}$ to $E$ three times, we get

$$\mathbf{u_3} = \bar{S}C\mathbf{u_2} + n^2 b\mathbf{i}, \ \ \mathbf{h_3} = \tilde{S}C\mathbf{h_2}.$$

Hence, after taking $k$ such that

$$\varepsilon_r < \frac{1}{n^k} \le n\varepsilon_r \ , \tag{3}$$

that is,

$$r > k \ge r-1 \tag{4}$$

and applying $W_{ij}$ to $E$ $k$ times ,we get $n^{2(k-1)}$ squares of side $1/n^k$ contained in $E_{ij}$ and

$$\Phi(\mathbf{h}_k\varepsilon_r^{-1}) - n^{2k} \le N(\varepsilon_r) \le \Phi(\mathbf{u}_k\varepsilon_r^{-1} + n^{2(k-1)}\mathbf{i})\left(\left[\frac{\varepsilon_r^{-1}}{n^k}\right]+1\right)^2 \tag{5}$$

, where

$$\mathbf{u}_k = \bar{S}C\mathbf{u}_{k-1} + n^{k-1}b\mathbf{i}, \ \ \mathbf{h}_k = \tilde{S}C\mathbf{h}_{k-1}.$$

Then

$$\mathbf{u}_k = (\bar{S}C)^{k-1}\mathbf{r} + (\bar{S}C)^{k-1}b\mathbf{i} + (\bar{S}C)^{k-2}nb\mathbf{i} + \cdots + (\bar{S}C)n^{k-2}b\mathbf{i} + n^{k-1}b\mathbf{i}, \tag{6}$$

$$\mathbf{h}_k = (\tilde{S}C)^{k-1}\mathbf{h}_1 \ . \tag{7}$$

Since $\tilde{S}C$, $\bar{S}C$ are non-negative irreducible matrix, from Frobenius' theorem there are strictly positive eigenvectors of $\tilde{S}C$ , $\bar{S}C$ 의 which corresponds to eigenvalues $\tilde{a} := \sum_{i,j=1}^{n} \tilde{s}_{ij}$ , $\bar{a} := \sum_{i,j=1}^{n} \bar{s}_{ij}$ , and we can choose strictly positive eigenvectors $\bar{\mathbf{a}}, \tilde{\mathbf{a}}$ which corresponds to eigenvalues $\bar{a}$ , $\tilde{a}$ so that

$$0 < \tilde{\mathbf{a}} < \mathbf{h_1}, \tag{8}$$

$$\mathbf{r} < \bar{\mathbf{a}}, \quad b\mathbf{i} < \bar{\mathbf{a}}n. \tag{9}$$

Then by (5),

$$N(\varepsilon_r) \le \Phi(\mathbf{u}_k \varepsilon_r^{-1} + n^{2(k-1)}\mathbf{i})\left(\left[\frac{\varepsilon_r^{-1}}{n^k}\right]+1\right)^2 \le \Phi(\mathbf{u}_k \varepsilon_r^{-1} + n^{2(k-1)}\mathbf{i})(n+1)^2$$

$$\le \Phi((\overline{S}C)^{k-1}\mathbf{r}\varepsilon_r^{-1} + (\overline{S}C)^{k-1}b\mathbf{i}\varepsilon_r^{-1} + (\overline{S}C)^{k-2}nb\mathbf{i}\varepsilon_r^{-1} + \cdots + (\overline{S}C)n^{k-2}b\mathbf{i}\varepsilon_r^{-1} +$$
$$+ n^{k-1}b\mathbf{i}\varepsilon_r^{-1} + n^{2(k-1)}\mathbf{i})(n+1)^2$$
$$\le \Phi((\overline{S}C)^{k-1}\overline{\mathbf{a}}\varepsilon_r^{-1} + (\overline{S}C)^{k-1}\overline{\mathbf{a}}\varepsilon_r^{-1}n + (\overline{S}C)^{k-2}\overline{\mathbf{a}}\varepsilon_r^{-1}n^2 + \cdots + (\overline{S}C)\overline{\mathbf{a}}\varepsilon_r^{-1}n^{k-1} +$$
$$+ \overline{\mathbf{a}}\varepsilon_r^{-1}n^k + n^{2(k-1)}\mathbf{i})(n+1)^2$$
$$\le (\overline{a}^{k-1}\Phi(\overline{\mathbf{a}})\varepsilon_r^{-1} + \overline{a}^{k-1}\Phi(\overline{\mathbf{a}})\varepsilon_r^{-1}n + \overline{a}^{k-2}\Phi(\overline{\mathbf{a}})\varepsilon_r^{-1}n^2 + \cdots + \overline{a}\Phi(\overline{\mathbf{a}})\varepsilon_r^{-1}n^{k-1} +$$
$$+ \Phi(\overline{\mathbf{a}})\varepsilon_r^{-1}n^k + n^{2(k-1)}\Phi(\mathbf{i}))(n+1)^2$$

$$\le (\overline{a}^{r-1}\overline{\mu}\varepsilon_r^{-1} + \overline{a}^{r-1}\overline{\mu}\varepsilon_r^{-1}n + \overline{a}^{r-2}\overline{\mu}\varepsilon_r^{-1}n^2 + \cdots + \overline{a}\overline{\mu}\varepsilon_r^{-1}n^{r-1} +$$
$$+ \overline{\mu}\varepsilon_r^{-1}n^r + \varepsilon_r^{-1}n^r)(n+1)^2$$

(10)

, where $\overline{\mu} = \Phi(\overline{\mathbf{a}})$.

On the other hands, since $0 < \tilde{s}_{ij} \le \overline{s}_{ij}$, $i, j = 1, \cdots, n$ , we have $\tilde{a} = \sum_{i,j=1}^{n} \tilde{s}_{ij} \le \sum_{i,j=1}^{n} \overline{s}_{ij} = \overline{a}$. If $\tilde{a} > n$, then $1 > \frac{n}{\tilde{a}} \ge \frac{n}{\overline{a}}$. Therefore,

$$N(\varepsilon_r) \le \overline{a}^{r-1}\varepsilon_r^{-1}\overline{\mu}\left(1 + n + \frac{n^2}{\overline{a}} + \cdots + \frac{n^r}{\overline{a}^{r-1}} + \frac{n^r}{\overline{\mu}\overline{a}^{r-1}}\right)(n+1)^2$$

$$= \overline{a}^{r-1}\varepsilon_r^{-1}\overline{\mu}\left(1 + \frac{n\left(1-\left(\frac{n}{\overline{a}}\right)^r\right)}{1-\frac{n}{\overline{a}}} + \frac{n^r}{\overline{\mu}\overline{a}^{r-1}}\right)(n+1)^2.$$

But since

$$\gamma := 1 + \frac{n\left(1-\left(\frac{n}{\overline{a}}\right)^r\right)}{1-\frac{n}{\overline{a}}} + \frac{n^r}{\overline{\mu}\overline{a}^{r-1}} > 0$$

, we have

$$\log N(\varepsilon_r) \le (r-1)\log\bar{a} + \log\varepsilon_r^{-1} + \log(\bar{\mu}\gamma(n+1)^2),$$

$$\frac{\log N(\varepsilon_r)}{-\log\varepsilon_r} \le (r-1)\frac{\log\bar{a}}{r\log n} + 1 + \frac{\log(\bar{\mu}\gamma(n+1)^2)}{-\log\varepsilon_r}$$

$$= \log_n^{\bar{a}} - \frac{\log\bar{a}}{r\log n} + 1 + \frac{\log(\bar{\mu}\gamma(n+1)^2)}{-\log\varepsilon_r}.$$

. where $\log x$ implies $\log_a^x$ with $a > 1$. Hence

$$\dim_B A = \lim_{\varepsilon_r \to 0}\frac{\log N(\varepsilon_r)}{-\log\varepsilon_r} \le 1 + \log_n^{\bar{a}} = 1 + \log_n^{\sum_{i,j=1}^{n}\bar{s}_{ij}}. \qquad (11)$$

By (5),

$$N(\varepsilon_r) \ge \Phi(\mathbf{h}_k\varepsilon_r^{-1}) - n^{2k} = \Phi((\tilde{S}C)^{k-1}\mathbf{h}_1\varepsilon_r^{-1}) - n^{2k}$$

$$\ge \Phi((\tilde{S}C)^{k-1}\tilde{\mathbf{a}}\varepsilon_r^{-1}) - n^{2k} = \tilde{a}^{k-1}\Phi(\tilde{\mathbf{a}})\varepsilon_r^{-1} - n^{2k}$$

$$\ge \tilde{a}^{r-2}\varepsilon_r^{-1}\tilde{\mu} - n^{2r} = \varepsilon_r^{-1}\tilde{a}^{r-2}\left(\tilde{\mu} - \frac{n^r}{\tilde{a}^{r-2}}\right),$$

Where $\tilde{\mu} = \Phi(\tilde{\mathbf{a}})$. But since $\tilde{a} > n$, there is $r_0$ such that for all $r(> r_0)$

$$\tilde{\mu} - \frac{n^r}{\tilde{a}^{r-2}} > 0.$$

Hence for $r(> r_0)$

$$\frac{\log N(\varepsilon_r)}{-\log\varepsilon_r} \ge 1 + (r-2)\frac{\log\tilde{a}}{r\log n} + \frac{\log\left(\tilde{\mu} - \frac{n^r}{\tilde{a}^{r-2}}\right)}{-\log\varepsilon_r}$$

$$\dim_B A = \lim_{\varepsilon_r \to 0}\frac{\log N(\varepsilon_r)}{-\log\varepsilon_r} \ge 1 + \log_n^{\tilde{a}} = 1 + \log_n^{\sum_{i,j=1}^{n}\tilde{s}_{ij}}. \qquad (12)$$

By (11), (12) if $\sum_{i,j=1}^{n}\tilde{s}_{ij} > n$, then we get

$$1 + \log_n^{\sum_{i,j=1}^{n}\tilde{s}_{ij}} \le \dim_B A \le 1 + \log_n^{\sum_{i,j=1}^{n}\bar{s}_{ij}}.$$

Proof of (2) . In the case of $\sum_{i,j=1}^{n}\bar{s}_{ij} = \bar{a} \le n$, by (10)

$$N(\varepsilon_r) \le (\bar{a}^{r-1}\bar{\mu}\varepsilon_r^{-1} + \bar{a}^{r-1}\bar{\mu}\varepsilon_r^{-1}n + \bar{a}^{r-2}\bar{\mu}\varepsilon_r^{-1}n^2 + \cdots + \bar{a}\bar{\mu}\varepsilon_r^{-1}n^{r-1} + + \bar{\mu}\varepsilon_r^{-1}n^r + \varepsilon_r^{-1}n^r)(n+1)^2$$

$$\le \varepsilon_r^{-1}n^r\bar{\mu}\left(n^{-1} + r + \frac{1}{\bar{\mu}}\right)(n+1)^2,$$

$$\frac{\log N(\varepsilon_r)}{-\log\varepsilon_r} \le 1 + \frac{\log n^r}{\log n^r} + \frac{\log\left(\bar{\mu}\left(n^{-1} + r + \frac{1}{\bar{\mu}}\right)(n+1)^2\right)}{-\log\varepsilon_r}$$

$$= 2 + \frac{1}{r}\log_n^{\left(\bar{\mu}\left(n^{-1}+r+\frac{1}{\bar{\mu}}\right)(n+1)^2\right)}$$

, therefore

$$\dim_B A = \lim_{\varepsilon_r \to 0}\frac{\log N(\varepsilon_r)}{-\log\varepsilon_r} \le 2.$$

On the other hands, since $A$ is the surface in $P^2$, we have $\dim_B A \ge 2$. Hence $\dim_B A = 2$.

(End of proof)

Remark1. The result of the theorem improves estimation of Box-counting dimension in [14]. In fact, let denote $\tilde{s} := \underset{E}{\mathrm{Min}}\,|s(x,y)|$, $\bar{s} := \underset{E}{\mathrm{Max}}\,|s(x,y)|$, and then since $\tilde{s} \le \tilde{s}_{ij} \le \bar{s}_{ij} \le \bar{s},\ (i, j = 1, \cdots, n)$, if $\tilde{s} > \frac{1}{n}$ then $\sum_{i,j=1}^{n} \underset{\sim}{s}_{ij} > n$,

$$1 + \log_n^{\sum_{i,j=1}^{n}\tilde{s}_{ij}} \ge 1 + \log_n^{\sum_{i,j=1}^{n}\tilde{s}} \ge 1 + \log_n^{n^2\tilde{s}} = 3 + \log_n^{\tilde{s}},$$

$$1 + \log_n^{\sum_{i,j=1}^{n}\bar{s}_{ij}} \le 1 + \log_n^{\sum_{i,j=1}^{n}\bar{s}} = 3 + \log_n^{\bar{s}}.$$

That is

$$3 + \log_n^{\tilde{s}} \le \dim_B A \le 3 + \log_n^{\bar{s}}.$$

And if $\tilde{s} \le \frac{1}{n}$, then since $\sum_{i,j=1}^{n}\bar{s}_{ij} \le n$, we have $\dim_B A = 2$. This is just estimation of Box-counting dimension in [14].

Remark2. for the surfaces that

$$f(x, y) = s_{ij}(x, y) f(L^{-1}(x, y)) + Q(x, y)$$

We get the same results as one of this theorem. In this case, $s_{ij}(x, y)$ is contraction function on $E_{ij}$ with $0 < |s_{ij}(x, y)| < 1$, and $L(x, y)$, $Q(x, y)$ are defined as in (2). In fact, we can choose $\bar{s}_{ij}$, $\tilde{s}_{ij}$ as $\bar{s}_{ij} = \underset{E_{ij}}{\mathrm{Max}} |s_{ij}(x, y)|$, $\underset{\sim}{s}_{ij} = \underset{E_{ij}}{\mathrm{Min}} |s_{ij}(x, y)|$ in the proof.